\documentclass[12pt]{article}
\usepackage{amssymb}
\usepackage{mathrsfs}
\usepackage{cite}
\usepackage{tikz}
\usepackage{amsmath}
\usepackage{bm}
\usepackage{graphicx}
\usepackage{pst-poly}  
\usepackage{pst-plot}  
\usepackage{multirow}
\usepackage{bigstrut}
\usepackage{url}  

\newtheorem{thm}{Theorem}[section]
\newtheorem{lem}[thm]{Lemma}
\newtheorem{Def}[thm]{Definition}
\newtheorem{cor}[thm]{Corollary}

\newenvironment{pf}[1][Proof]{\noindent\textbf{#1.} }{\hfill\rule{1mm}{2mm}}

\makeatletter \@addtoreset{equation}{section} \makeatother

\begin{document}
\title{\bf On conditional fault tolerance of hierarchical cubic networks\thanks{The work was supported by NNSF of China (11571044, 11601041),
Young Talent Fund of Yangtze University (2015cqr23).}}

\author
{Xiang-Jun Li$^a$ \quad Min Liu$^a$\quad Zheng Yan$^a$ \quad
Jun-Ming Xu$^b$ \footnote{Corresponding author: xujm@ustc.edu.cn
(J.-M. Xu) }\\
\\
 {\small $^a$School of Information and Mathematics,}\\
 {\small Yangtze University, Jingzhou, Hubei, 434023, China}\\
 {\small $^b$School of Mathematical Sciences}\\
 {\small University of Science and Technology of China, Hefei, 230026, China}}
\date{}
\maketitle

\begin{abstract}

This paper considers the conditional fault tolerance, $h$-super
connectivity $\kappa^{h}$ and $h$-super edge-connectivity
$\lambda^{h}$ of the hierarchical cubic network $HCN_n$, an
attractive alternative network to the hypercube, and shows
$\kappa^h(HCN_n)=\lambda^h(HCN_n)=2^h(n+1-h)$ for any $h$ with
$0\leq h\leq n-1$. The results imply that at least $2^h(n+1-h)$
vertices or edges have to be removed from $HCN_n$  to make it
disconnected with no vertices of degree less than $h$, and
generalize some known results.

\vskip6pt


\vskip6pt


\end{abstract}

\section{Introduction}
\quad \ It is well known that interconnection networks play an important
role in parallel computing/communication systems.
An interconnection network can be modeled by a graph $G=(V, E)$, where $V$ is the set
of processors and $E$ is the set of communication links in the
network.

The $n$-dimensional hypercube $Q_n$ is a graph whose vertex-set
consists of all binary vectors of length $n$, with two vertices
being adjacent whenever the corresponding vectors differ in exactly
one coordinate. For its regularity, symmetry, high connectivity,
logarithmic diameter and simple routing, the hypercube becomes one
of the most popular, versatile and efficient topological structures
of interconnection networks~\cite{l92}.

However, the hypercube has been considered unsuitable for building
large systems since the relatively high vertex-degree results in an
additional difficulty in interconnection. To make up for these
defects, as an alternative to the hypercube network, many variations
of the hypercube network are proposed in the literature. One of them
is the hierarchical cubic networks $HCN_n$ proposed by Ghose and
Desai~\cite{GD95}, which is feasible to be implemented with
thousands of or more processors, while retaining a good properties
from the hypercubes, such as regularity, symmetry and logarithmic
diameter. Compared with the hypercube of the same size, the
hierarchical cubic network requires only about half the number of
edges and provides a lower diameter~\cite{Chc96,fcd02,GD95,yp98}.

In real networks, since the fault of vertices and edges are
inevitable, measuring the fault tolerance in networks are very
important. The traditional connectivity is a good measurement for
the fault tolerance of networks. The {\it connectivity} $\kappa(G)$
(resp. {\it edge-connectivity} $\lambda(G)$ ) of $G$ is defined as
the minimum number of vertices (resp. edges) whose removal from $G$
results in a disconnected graph. The connectivity $\kappa(G)$ and
edge-connectivity $\lambda(G)$ of a graph $G$ are two important
measurements for fault tolerance of the network since the larger
$\kappa(G)$ or $\lambda(G)$ is, the more reliable the network is
(see~\cite{x01}).

However, the definitions of $\kappa(G)$ and $\lambda(G)$ are
implicitly assumed that any subset of system components is equally
likely to be faulty simultaneously, which may not be true in real
applications, thus they underestimate the reliability of the
network. To overcome such a shortcoming, Harary~\cite{h83}
introduced the concept of conditional connectivity by appending some
requirements on connected components, Latifi {\it et
al.}~\cite{lhm94} specified requirements and proposed the concept of
the restricted $h$-connectivity. These parameters can measure fault
tolerance of an interconnection network more accurately than the
classical connectivity. The concepts stated here are slightly
different from theirs (see~\cite{xlmh05}).

For a graph $G$, $\delta(G)$ denotes its minimum vertex-degree. A
subset $S\subset V(G)$ (resp. $F\subset E(G)$) is called an {\it
$h$-vertex-cut} (resp. {\it edge-cut}), if $G-S$ (resp. $G-F$) is
disconnected and $\delta(G-S)\geq h$. The {\it $h$-super
connectivity} $\kappa^{h}(G)$ (resp. {\it $h$-super
edge-connectivity} $\lambda^{h}(G)$) of $G$ is defined as the
cardinality of a minimum $h$-vertex-cut (resp. $h$-edge-cut) of $G$.
It is clear that $\kappa^{0}(G)=\kappa(G)$ and
$\lambda^{0}(G)=\lambda(G)$.

For an arbitrarily given graph $G$ and any integer $h$, determining
the exact values of $\kappa^{h}(G)$ and $\lambda^{h}(G)$ is quite
difficult, no polynomial algorithm to compute them has been yet
known so far. In fact, the existence of $\kappa^{h}(G)$ and
$\lambda^{h}(G)$ is an open problem for a general graph $G$ and
$h\geq 1$. The main interest of the researchers is to determine the
values of $\kappa^{h}$ and $\lambda^{h}$ for some well-known classes
of networks and any $h$. For a long time, almost all of the research
on this topics has been focused on some small $h$'s, only the
hypercube network, its $\kappa^{h}$ and $\lambda^{h}$ were
determined~\cite{oc93,wg98,x00c} for any $h$ with $0\leq h\leq n-2$.

In recent years, some new methods and techniques have been
discovered, from which $\kappa^{h}$ and $\lambda^{h}$ have been
determined for some well-known classes of networks and for any $h$.
For example, $\kappa^{h}$ and $\lambda^{h}$ were determined for star
networks~\cite{lx14}, $(n, k)$-star networks~\cite{lx12,lx14,
lgyx17} and exchanged hypercubes~\cite{lx13}; $\kappa^{h}$ was
determined for $(n, k)$-arrangement networks~\cite{lz12}, exchanged
crossed cubes~\cite{n17} and locally twisted cubes~\cite{wh17};
$\lambda^{h}$ was determined for hypercube-like
networks~\cite{lx12c}.

This paper is interested in the hierarchical cubic network $HCN_n$.
Chiang and Chen~\cite{Chc96} determined
$\kappa(HCN_n)=\lambda(HCN_n)=n+1$, Zhou {\it et al.}~\cite{z16}
proved that $\kappa^1(HCN_n)=2n$ and $\kappa^2(HCN_n)=4(n-1)$. We
generalize these results by proving that
$\kappa^h(HCN_n)=2^h(n+1-h)$ for any $h$ with $0\leq h\leq n-1$, and
$\lambda^h(HCN_n)=2^h(n+1-h)$ for any $h$ with $0\leq h\leq n$.

The rest of the paper is organized as follows. In Section 2, we
recall the structure of $HCN_n$ and some lemmas used in our
proofs. The main proof of the result is in Section 3.
Conclusions are in Section 4.

For graph terminology and notation not defined here we follow
Xu~\cite{x01}. For a subset $X$ of vertices in $G$, we do not
distinguish $X$ and the induced subgraph $G[X]$.


\section{Definitions and lemmas}

Let $V_n$ be the set of binary sequence of length $n$, i.e.,
$V_n=\{x_1x_2\cdots x_n:\ x_i\in \{0,1\}, 1\leq i\leq n\}$. For
$x=x_1x_2\cdots x_n\in V_n$, the element $\bar x=\bar x_1\bar
x_2\cdots \bar x_n\in V_n$ is called the bitwise complement of $x$,
where $\bar x_i=\{0,1\}\setminus\{x_i\}$ for each
$i\in\{1,2,\ldots,n\}$.

A hypercube network $Q_n$ is an $n$-dimensional cube, shortly
$n$-cube, its vertex-set $V_n$, and two vertices being linked by an
edge if and only if they differ exactly in one coordinate. For the
sake of simplicity, we use $xQ_n$ to denote the Cartesian product
$\{x\}\times Q_n$ of a vertex $x$ and a hypercube network $Q_n$.

\begin{Def}\textnormal{(\cite{GD95})}\label{def2.1}
An $n$-dimensional hierarchical cubic network $HCN_n$ with
vertex-set $V_n\times V_n$ is obtained from $2^n$ $n$-cubes $\{x
Q_n: x\in V_n\}$ by adding edges between two $n$-cubes, called {\it
crossing edges}, according to the following rule: A vertex $(x, y)$
in $xQ_n$ is linked to

(1) $(y, x)$ in $yQ_n$ if $x \neq y$ or

(2) $(\bar{x}, \bar{y})$ in $\bar xQ_n$ if $x = y$.

\noindent The vertex $(y, x)$ in $yQ_n$ or $(\bar{x}, \bar{y})$ in
$\bar xQ_n$ is called an external neighbor of $(x, y)$ in $xQ_n$.
\end{Def}

\begin{figure}[ht]
\begin{center}
\psset{unit=25pt}
\begin{pspicture}(-4,-3.5)(4,4)

\cnode(-3,3){3pt}{a1}\rput(-3.7,3.4){\scriptsize(01,01)}
\cnode(-1,3){3pt}{a2} \rput(-1.7,3.4){\scriptsize(01,11)}
\cnode(-1,1){3pt}{a3} \rput(-1.7,1.3){\scriptsize(01,10)}
\cnode(-3,1){3pt}{a4} \rput(-3.7,1.3){\scriptsize(01,00)}
\ncline{a1}{a2}\ncline{a2}{a3}\ncline{a3}{a4}\ncline{a1}{a4}

\cnode(1,3){3pt}{b1} \rput(1.4,3.4){\scriptsize(11,01)}
\cnode(3,3){3pt}{b2} \rput(3.8, 3.4){\scriptsize(11,11)}
\cnode(3,1){3pt}{b3}\rput(3.8,1.3){\scriptsize(11,10)}
\cnode(1,1){3pt}{b4} \rput(0.3,1.3){\scriptsize(11,00)}
\ncline{b1}{b2}\ncline{b2}{b3}\ncline{b3}{b4}\ncline{b1}{b4}

\cnode(-3,-1){3pt}{c1}  \rput(-3.7,-.7){\scriptsize(00,01)}
\cnode(-1,-1){3pt}{c2}  \rput(-1.7,-.7){\scriptsize(00,11)}
\cnode(-1,-3){3pt}{c3}  \rput(-1.7,-2.7){\scriptsize(00,10)}
\cnode(-3,-3){3pt}{c4}  \rput(-3.7,-2.7){\scriptsize(00,00)}
\ncline{c1}{c2}\ncline{c2}{c3}\ncline{c3}{c4}\ncline{c1}{c4}

\cnode(1,-1){3pt}{d1}  \rput(0.3,-1.3){\scriptsize(10,01)}
\cnode(3,-1){3pt}{d2} \rput(3.8,-0.7){\scriptsize(10,11)}
\cnode(3,-3){3pt}{d3} \rput(3.8,-2.7){\scriptsize(10,10)}
\cnode(1,-3){3pt}{d4}\rput(1.8,-2.7){\scriptsize(10,00)}
\ncline{d1}{d2}\ncline{d2}{d3}\ncline{d3}{d4}\ncline{d1}{d4}

\ncline[linecolor=red]{a4}{c1} \ncline[linecolor=red]{a3}{d1}
\ncline[linecolor=red]{a2}{b1} \ncline[linecolor=red]{c3}{d4}
\ncline[linecolor=red]{b4}{c2} \ncline[linecolor=red]{b3}{d2}

\nccurve[angleA=-30,angleB=120,linecolor=red]{a1}{d3}
\nccurve[angleA=30,angleB=-120,linecolor=red]{c4}{b2}

\end{pspicture}
\caption{\label{f1}\footnotesize {$2$-dimensional hierarchical cubic
network $HCN_2$}}
\end{center}
\end{figure}
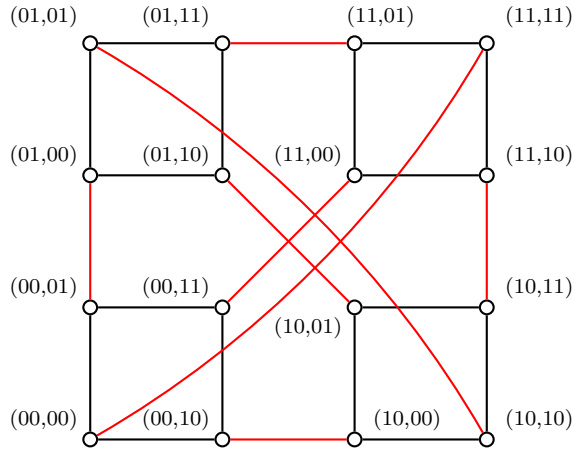

A $2$-dimensional hierarchical cubic network $HCN_2$ is shown in
Fig.~\ref{f1}, where the red edges are the crossing edges in
$HCN_2$.

Clearly, $HCN_n$ is an $(n+1)$-regular graph. Chiang and
Chen~\cite{Chc96} determined its connectivity and edge-connectivity.

\begin{lem}\label{lem2.2}\textnormal{(\cite{Chc96})}
$\kappa(HCN_n)=\lambda(HCN_n)=n+1$.
\end{lem}

From Definition~\ref{def2.1}, it is easy to obtain the following
property about crossing edges in $HCN_n$.

\begin{lem}\label{lem2.3}
(1) There are two crossing edges between two $n$-cubes $x Q_n$ and
$y Q_n$ if and only if $x$ and $y$ are complementary; otherwise
there is only one crossing edge. (2) The set of crossing edges
consists of a perfect matching of $HCN_n$.
\end{lem}

Since $HCN_n$ is made up of $2^n$ $n$-cubes and a perfect matching,
some properties on an $n$-cube $Q_n$ are very useful for the proofs
of our main results.

\begin{lem}\label{lem2.4}\textnormal{(\cite{oc93,wg98,x00c}  )}
$\kappa^h(Q_n)=2^h(n-h)$  for any $h$ with $0\leq h\leq n-2$, and
$\lambda^h(Q_n)=2^h(n-h)$ for any $h$ with $0\leq h\leq n-1$.
\end{lem}

\begin{lem}\label{lem2.5}\textnormal{(\cite{wg98})}
If $X$ is a subgraph in $Q_n$ and $\delta(X)\geq h$,
 then $|X|\geq 2^h$.
\end{lem}

For a subgraph $X$ in $Q_n$, $N_n(X)$ denotes the set of neighbors
of $X$ in $Q_n - X$.

\begin{lem}\label{lem2.6}
If $X$ is a subgraph in $Q_n$ and $\delta(X)\geq h$, then
$|X|+|N_n(X)|\geq 2^h(n-h)$ for any $h$ with $0\leq h\leq n-1$ and
$n\ge 1$.
\end{lem}

\begin{pf}
For $n=1$, $Q_1\cong K_2$, the conclusion holds clearly. Assume
$n\ge 2$ below. The proof proceeds by induction on $h\geq 0$ for a
fixed $n$. Since $Q_n$ is $n$-regular, for any non-empty subgraph
$X$ of $Q_n$, $|X|+|N_n(X)|\geq n+1$, and so the conclusion is true
for $h=0$. Assume the induction hypothesis for $h-1$ with $h \geq
1$.

It is well known that $Q_n$ can be expressed as $Q_n=L\odot_i R$,
where $L$ and $R$ are two $(n-1)$-cubes induced by the vertices with
$i$-th coordinate is 0 and $1$, respectively, the set of edges
between $L$ and $R$ consists of a perfect matching in $Q_n$ (see
Xu~\cite{x01}).

Let $X$ be a subgraph in $Q_n$ with $\delta(X)\geq h$. Then
$E(X)\ne\emptyset$ since $h\geq 1$. Arbitrarily take an edge $e$ of
$X$, and assume that two end-vertices of $e$ differ in only the
$i$-th coordinate. Let $Q_n=L\odot_i R$. Then $X\cap L\ne\emptyset$
and $X\cap R\ne\emptyset$.

Let $X_0=X\cap L, X_1=X\cap R$. Since $\delta(X)\geq h$ in $Q_n$ and
the set of edges between $L$ and $R$ is a matching, $\delta(X_0)\geq
h-1$ in $L$ and $\delta(X_1)\geq h-1$ in $R$. Using the induction
hypothesis in $L$ and $R$ respectively, we have
 $$
\begin{array}{l}
|X_i|+|N_{n-1}(X_i)| \geq 2^{h-1}(n-h)\ \ {\rm for\ each}\ i\in
\{0,1\}.
\end{array}
 $$
It follows that
  $$
\begin{array}{l}
|X|+|N_{n}(X)|\geq |X_0|+|N_{n-1}(X_0)|+|X_1|+|N_{n-1}(X_1)|\geq
2^{h}(n-h).
\end{array}
 $$

By the induction principle, the lemma follows.
\end{pf}

\section{Main results}

\begin{lem}\label{lem2.1}
For $n\ge 1$, $\kappa^h(HCN_n)\leq 2^h(n+1-h)$ for any $h$ with
$0\leq h\leq n-1$, and  $\lambda^h(HCN_n)\leq 2^h(n+1-h)$ for any
$h$ with $0\leq h\leq n$.
\end{lem}
\begin{pf}
For $n=1$, $HCN_1\cong C_4$, a cycle of length 4, the conclusion
holds clearly. Assume $n\ge 2$ below. Let $x_1Q_n, x_2Q_n,\ldots,
x_{2^n}Q_n$ be $2^n$ $n$-cubes in $HCN_n$. For a fixed $h$ with
$0\leq h\leq n-1$, let $x_1Q_{h}$ be a subgraph in $x_1Q_n$ induced
by the vertices with the rightmost $(n-h)$ bits $0$s of the second
component, $S$ be the neighbors of $x_1Q_{h}$ in $HCN_n - x_1Q_{h}$.
Then $HCN_n- S$ is disconnected.

On the one hand, by the choice of $Q_h$, $S$ must contain all
vertices with exactly one $1$ in the rightmost $(n-h)$ coordinates
of the second component, such vertices have exactly $2^h(n-h)$. On
the other hand, $S$ must contain all external neighbors of $x_1Q_h$,
such external neighbors have exactly $2^h$. Thus,
$|S|=2^h(n-h)+2^h=2^h(n+1-h)$.

\vskip6pt

We now need to prove that $S$ is an $h$-vertex-cut, i.e., each
vertex in $HCN_n-S$ has at least $h$ neighbors.

We first show that $|S\cap V(x_jQ_n)|\leq 1$ for each $j\ne 1$. On
the contrary, suppose that $|S\cap V(x_{j_0}Q_n)|=2$ for some
$j_0\ne 1$. Then there are two crossing edges, say $e_1$ and $e_2$,
between $x_1Q_h$ and $x_{j_0}Q_n$, and so $j_0=\bar x_1$ by
Lemma~\ref{lem2.3}. By Definition~\ref{def2.1}, two of end-vertices
of $\{e_1,e_2\}$ in $x_1Q_h$ is certainly $(x_1,x_1)$ and $(x_1,\bar
x_1)$. Since the distance between $(x_1,x_1)$ and $(x_1,\bar x_1)$
is $n$, we have $n\leq h$, a contradiction. It follows that $|S\cap
V(x_jQ_n)|\leq 1$ for each $j\ne 1$.

\vskip6pt

For any $j\ne 1$, let $z$ be a vertex in $x_jQ_n-S$. Then $z$ has at
most one neighbor in $S\cap V(x_jQ_n)$ since $|S\cap V(x_jQ_n)|\leq
1$. By Definition~\ref{def2.1}, $z$ has at most one neighbor in
$HCN_n-x_{j}Q_n$. Thus, the vertex $z$ has at least $(n+1)-2\geq h$
neighbors in $HCN_n-S$.

Let $S_1=V(x_1Q_n)\cap S$ and $T_1=V(x_1(Q_n-Q_h)-S_1)$. All that's
left is to prove that each vertex in $x_1Q_n-S_1$ has at least $h$
neighbors in $HCN_n - S$. It is clear that each vertex in $x_1Q_h$
has $h$ neighbors in $HCN_n-S$ by the choice of $x_1Q_h$.

If $T_1$ is empty then we have done. Assume $T_1\ne\emptyset$ and
let $w\in T_1$. Then $h\le n-2$. If $w$ has no neighbors in $S_1$,
then it has at least $n$ neighbors in $HCN_n - S$. Suppose that $w$
has neighbors in $S_1$. By the choice of $x_1Q_h$, there is exactly
one $1$ in the rightmost $(n-h)$ coordinates of the second component
of each vertex in $S_1$, and so there are exactly two $1$s in the
rightmost $(n-h)$  coordinates of the second component of $w$, which
implies that $w$ has at most two neighbors in $S_1$. Thus $w$ has at
least $(n-2)\geq h$ neighbors in $HCN_n - S$.

\vskip6pt

From the above discussions, each vertex of $HCN_n - S$ has at least
$h$ neighbors within. Therefore, $S$ is an $h$-vertex-cut in
$HCN_n$, and so $\kappa^h(HCN_n)\leq |S|=2^h(n+1-h)$.

\vskip6pt

Let $F$ be the set of edges between $x_1Q_{h}$ and $S$. Then $HCN_n-
F$ is disconnected. From the above discussions, it is easy to see
that $F$ is an $h$-edge-cut in $HCN_n$ and $|F|=|S|$. Thus
$\lambda^h(HCN_n)\leq |F|=|S|=2^h(n+1-h)$.

The lemma follows.
\end{pf}

 \begin{thm}\label{thm3.2}
For $n\ge 1$, $\kappa^h(HCN_n)=2^h(n+1-h)$ for any $h$ with $0\leq
h\leq n-1$, and  $\lambda^h(HCN_n)=2^h(n+1-h)$ for any $h$ with
$0\leq h\leq n$.
\end{thm}
\begin{pf}
For $n=1$, $HCN_1\cong C_4$, a cycle of length 4, the conclusion
holds clearly. Assume $n\ge 2$ below. By Lemma~\ref{lem2.2} and
Lemma~\ref{lem2.1}, we only to show $\kappa^h(HCN_n)\geq 2^h(n+1-h)$
for any $h$ with $1\leq h\leq n-1$, and  $\lambda^h(HCN_n)\geq
2^h(n+1-h)$ for any $h$ with $1\leq h\leq n$.

To the end, let $F$ be a minimum $h$-vertex-cut (or $h$-edge-cut) of
$HCN_n$, $X$ be the vertex-set of the minimum connected component of
$HCN_n-F$, and let
 $$
 Y=\left\{\begin{array}{ll}
 V(HCN_n-X\cup F)\ & \text{if $F$ is a vertex-cut};\\
 V(HCN_n-X) \ & \text{if $F$ is an edge-cut}.
 \end{array}\right.
 $$

Let $H_1, H_2,\ldots, H_{2^n}$ be $2^n$ $n$-cubes in $HCN_n$. For
any $i\in \{1,2,\ldots,2^n\}$, let
 \begin{displaymath}{}
 \begin{array}{ll}
 X_i=X\cap V(H_i), Y_i=Y\cap V(H_i), \\
 F_i=\left\{
 \begin{array}{rl}
 F\cap V(H_i) &\ \text{if  $F$ is a vertex-cut};\\
  F\cap E(H_i) &\ \text{if $F$ is an edge-cut}.
 \end{array}\right.\\
 F_C=\left\{
 \begin{array}{rl}
   {\emptyset \quad \ \ \ \ } &\ \text{if  $F$ is a vertex-cut};\\
  F\cap \left(\bigcup_{i\ne j}E(H_i,H_j)\right) &\ \text{if $F$ is an edge-cut}.
 \end{array}\right.
   \end{array}
\end{displaymath}
where $E(H_i,H_j)$ denotes the set of edges between $H_i$ and $H_j$
for $i\ne j$. Let
\begin{displaymath}
\begin{array}{l}
 J_X=\{i\in \{1,2,\ldots,2^n\}:\ X_i\ne\emptyset\},\\
 J_Y=\{i\in \{1,2,\ldots,2^n\}:\ Y_i\not=\emptyset\}\ \ {\rm and} \\
 J_0 =J_X\cap J_Y.
 \end{array}
 \end{displaymath}

Clearly, if $J_0\ne\emptyset$ then $X_i\ne \emptyset$ and $Y_i\ne
\emptyset$ for each $i\in J_0$. By the choice of $F$, every vertex
in $X_i\cup Y_i$ has at least $h$ neighbors in $HCN_n-F$, at most
one of them is an external neighbor. This fact implies that $F_i$ is
an $(h-1)$-vertex-cut of $H_i$ if $F$ is a vertex-cut, or an
$(h-1)$-edge-cut of $H_i$ if $F$ is an edge-cut. Since $H_i$ is an
$n$-cube and $h-1\geq 0$, by Lemma~\ref{lem2.4} we have
  \begin{equation}\label{e3.1}
  |F_i|\geq 2^{h-1}(n+1-h) \ {\rm \ for\ each}\ i\in
  J_0,
 \end{equation}
and by Lemma~\ref{lem2.5} we have
 \begin{equation}\label{e3.2}
 |X_i|\geq 2^{h-1}\ \  {\rm and}\ \ |Y_i|\geq 2^{h-1}  \ {\rm \ for\ each}\ i\in
  J_0.
 \end{equation}

If $h=n$ then $F$ is an $n$-edge-cut. We will prove $|F|\geq 2^n$.

If $J_0=\emptyset$, then $F$ is only consists of crossing edges. Let
$G$ be a contracting graph of $HCN_n$, obtained by contracting each
$n$-cube $H_i$ in $HCN_n$ as a single vertex $x_i$ and by removing
all loops. It is easy to see that $G$ is a complete graph $K_{2^n}$
plus a perfect matching, and $F$ is an edge-cut of $G$. Thus,
$|F|\geq\lambda(G)=2^n$.

If $J_0\ne\emptyset$ then, $|F_i|\geq 2^{n-1}$ for $i\in J_0$ by
(\ref{e3.1}). Combining (\ref{e3.2}) with $H_i\cong Q_{n}$, we have
$|F_i|=|X_i|=|Y_i|=2^{n-1}$ and $X_i$ is $(n-1)$-regular for each
$i\in J_0$. Without loss of generality, assume $1\in J_0$. Since
$\delta(X)\ge n\geq 2$ and $X_1$ is $(n-1)$-regular, all external
neighbors of $X_1$ are certainly in $X\setminus X_1$. So $|J_X|\geq
|X_1|+1=2^{n-1}+1$. Since $|X_i|=|Y_i|$ for each $i\in J_0$, by the
minimality of $X$, we have $|J_Y|\geq |J_X|\geq 2^{n-1}+1$. Since
$|J_X \cup J_Y|=2^n$, we have $|J_0|=|J_X|+|J_Y|-|J_X \cup J_Y|\geq
2$. Thus, $|F|\geq \sum_{i\in J_0}|F_i|\geq 2\times 2^{n-1}=2^n$.

\vskip6pt

In the following discussion, we assume $1\leq h\leq n-1$ and need to
show that
 \begin{equation}\label{e3.3}
  |F| \geq 2^h(n+1-h) \ {\rm \ for}\ 1 \leq h \leq n-1.
 \end{equation}

If $|J_0| \geq 2$ then, by (\ref{e3.1}), we have that
  $$
\begin{array}{rl}
 |F| &\geq  \sum\limits_{i\in J_0}|F_i|\geq 2\times2^{h-1}(n+1-h)= 2^h(n+1-h).
\end{array}
$$
Thus, (\ref{e3.3}) holds if $|J_0| \geq 2$. Assume $0\leq |J_0| \leq
1$ below.

Let $a=|J_X \setminus J_0|, b=|J_Y \setminus J_0|,
c=|\{1,\ldots,2^n\}\setminus(J_X\cup J_Y)|$. By the choice of $X$
with minimum cardinality, we have $a\le b$. If $c\geq 1$, then there
exists some $i$ such that $V(H_i) \subseteq F $ and $F$ is a
vertex-cut, therefore $|F|\geq 2^n \geq 2^h(n+1-h)$ for $h\leq n-1$,
and so (\ref{e3.3}) holds. Next, assume $c=0$, that is,
$a+b+|J_0|=2^n$.

If $a\geq 1$ and $b\geq 1$ then, by Lemma~\ref{lem2.3}, for $j_1\in
J_X\setminus J_0,j_2\in J_Y\setminus J_0$, there is at least one
crossing edge between $H_{j_1}$ and $H_{j_2}$, and so there are at
least $ab$ crossing edges between $\cup_{j_1\in J_X\setminus
J_0}H_{j_1}$ and $\cup_{j_2\in J_Y\setminus J_0}H_{j_2}$. Each of
these crossing edges must be in $F$ if $F$ is an edge-cut, or one of
its end-vertices must be in $F$ if $F$ is a vertex-cut. Therefore,
we have

 \begin{equation}\label{e3.4}
  \sum_{i\in J_X\cup J_Y\setminus J_0} |F_i|+|F_C|\geq
  \sum_{i\in J_X\setminus J_0, j\in J_Y\setminus J_0}|E({H_i, H_j})|\geq ab.
 \end{equation}

We consider two cases depending on $|J_0|=0$ or $|J_0|=1$.

{\bf Case 1.}\   $|J_0| =0$.

In this case, $a\geq 1$. If $a\geq 2$, by (\ref{e3.4}) we have
 $$
 |F|\geq \sum_{i\in J_X\cup J_Y} |F_i|+|F_C| \geq ab = a(2^n-a)\geq 2^n\geq 2^h(n+1-h).
$$

If $a=1$, without loss of generality assume $J_X=\{1\}$, then
$X_1\subseteq V(H_1)$ if $F$ is a vertex-cut or $X_1=V(H_1)$ if $F$
is an edge-cut. If $F$ is a vertex-cut, then all external neighbors
of $X_1$ and all vertices in $V(H_1-X_1)$ are contained in $F$, and
so $|F|\geq |V(H_1)|=2^n$. If $F$ is an edge-cut, then all crossing
edges incident with $H_1$ are contained in $F$, and so $|F| \geq
|V(H_1)|=2^n$. Whether $F$ is a vertex-cut or an edge-cut, we have
$|F| \geq 2^n\geq 2^h(n+1-h)$ for $1\leq h\leq n-1$.

 \vskip6pt

{\bf Case 2.}\   $|J_0| =1$.

In this case, $a\geq 0$ and $b=2^n-a-1$. Without loss of generality,
we assume $J_0=\{1\}$.

 If $a\geq 1$, combining (\ref{e3.1}) and (\ref{e3.4}), we have

 $$
\begin{array}{rl}
 |F|&\geq |F_1|+\sum_{i\in J_X\cup J_Y\setminus J_0} |F_i|+|F_C| \\
 &\geq  2^{h-1}(n+1-h)+a(2^n-a-1)\\
 &\geq 2^{h-1}(n+1-h)+2^{n-1} \\
 &\geq 2^{h-1}(n+1-h)+2^{h-1}(n+1-h)\\
 &\geq 2^{h}(n+1-h).
\end{array}
$$

If $a=0$, then $J_X=J_0=\{1\}$. Since $\delta(X)\geq h$ and $H_1$ is
an $n$-cube, by Lemma~\ref{lem2.6} $|X|+|N_{H_1}(X)| \geq
2^{h}(n+1-h)$. If $F$ is a vertex-cut, then $N_{HCN_n}(X)\subset F$,
and so
 $$
|F|\geq |N_{HCN_n}(X)|\geq |X|+|N_{H_1}(X)|\geq 2^{h}(n+1-h).
 $$
If $F$ is an edge-cut then $F_1$ is the set of edges between $X$ and
$N_{H_1}(X)$, and so $|F_1|\geq |N_{H_1}(X)|$. Note that $|F_C|\geq
|X|$ since $a=0$. It follows that
 $$
|F|\geq |F_C|+|F_1|\geq |X|+|N_{H_1}(X)|\geq 2^{h}(n+1-h).
 $$

The theorem follows.
\end{pf}

\vskip6pt

Zhou {\it et al.}\cite{z16} determined $\kappa^1(HCN_n)$ and
$\kappa^2(HCN_n)$, which can be obtained from Theorem~\ref{thm3.2}
by setting  $h = 1,2$ respectively.

\begin{cor} \textnormal{(Zhou {\it et al.} \cite{z16} )}
 $\kappa^1(HCN_n)=2n$ and $\kappa^2(HCN_n)=4(n-1)$ for  $n\geq 3$.
\end{cor}

\section{Conclusions}
In this paper, we investigate the refined measure, $k$-super
connectivity $\kappa^{h}$ and $k$-super edge-connectivity
$\lambda^{h}$ for the fault tolerance of a network. For the
hierarchical cubic network $HCN_n$, which is an attractive
alternative network to the hypercube, we prove
$\kappa^h(HCN_n)=2^h(n+1-h)$ for any $h$ with $0\leq h\leq n-1$, and
$\lambda^h(HCN_n)=2^h(n+1-h)$ for any $h$ with $0\leq h\leq n$,
which implies that at least $2^h(n+1-h)$ vertices or edges have to
be removed from $HCN_n$ to make it disconnected with no vertices of
degree less than $h$. When the hierarchical cubic networks $HCN_n$
is used to model the topological structure of a large-scale parallel
processing system, these results can provide a more accurate measure
for the fault tolerance of the system. 

\end{document}